\documentclass[12pt,a4paper]{amsart}
 \usepackage{amssymb}
 \usepackage{multirow}  
 \usepackage{hhline}
 \usepackage{ifthen}
 \usepackage{amscd}
 \usepackage{graphicx}

 \theoremstyle{plain}
  \newtheorem{thm}{Theorem}[section]
  \newtheorem{lem}[thm]{Lemma}
  \newtheorem{slem}[thm]{Sublemma}
  \newtheorem{cor}[thm]{Corollary}

  \newtheorem{ass}[thm]{Assertion}

\theoremstyle{definition}
  \newtheorem{defn}[thm]{Definition}

\theoremstyle{remark}
  \newtheorem{rem}[thm]{Remark}

  \newtheorem{ack}{Acknowledgment}

\numberwithin{equation}{section}

\newenvironment{proclaim}[1]
  {\par \vspace{1.5ex} \noindent #1 \ }
  {\par \vspace{1.5ex}}

\DeclareMathOperator{\vol}{vol}                 
\DeclareMathOperator{\area}{area}               
\DeclareMathOperator{\diam}{diam}               
\DeclareMathOperator{\interior}{int}            
\DeclareMathOperator{\strrad}{str.rad}          

\newcommand{\field}[1]{\mathbb{#1}}

\newcommand{\R}{\field{R}}                      

\newcommand{\const}{\text{\rm const}}           

\DeclareMathOperator{\CCap}{Cap}

\newcommand{\cal}{\mathcal}

\begin{document}


\title{Volume collapsed three-manifolds\\
   with a lower curvature bound
}

\author{Takashi Shioya}
\author{Takao Yamaguchi}

\address{Mathematical Institute, Tohoku University, Sendai 980-8578, JAPAN}

\address{Institute of Mathematics, University of Tsukuba, Tsukuba
  305-8571,\newline
\hspace*{0.4cm} JAPAN}

\email[T. Shioya]{shioya@math.tohoku.ac.jp}
\email[T. Yamaguchi]{takao@math.tsukuba.ac.jp}

\date{}

\subjclass[2000]{Primary 53C20, 53C23; Secondary 57N10, 57M99}

\keywords{the Gromov-Hausdorff convergence, Alexandrov spaces, topology
  of three-manifolds, graph manifolds}

\begin{abstract}
  In this paper we determine the topology
  of three-dimensional closed orientable Riemannian manifolds 
   with a uniform lower bound
  of sectional curvature whose volume is sufficiently small.
\end{abstract}

\maketitle


\section{Introduction} \label{sec:intro}

As a continuation of our investigation \cite{SyYm:3mfd}
of collapsing three-manifolds
with  a lower curvature bound and an upper diameter bound, 
we study the topology of a three-dimensional closed Riemannian 
manifold with a lower curvature bound whose volume is 
sufficiently small, where we assume no upper diameter bound.

A closed three-manifold is called a {\em graph manifold}
if it is a finite gluing of Seifert fibered spaces
along their boundary tori.

\begin{thm} \label{thm:graph}
There exist small positive numbers $\epsilon_0$ and $\delta_0$ such that 
if a closed orientable three-manifold $M$ has a Riemannian 
metric with sectional curvature $K\ge -1$ and $\vol(M)< \epsilon_0$, 
then one of the following holds$:$
 \begin{enumerate}
  \item $M$ is homeomorphic to a graph manifold$;$
  \item $\diam(M) < \delta_0$ and $M$ has finite fundamental 
     group.
 \end{enumerate}
\end{thm}

It was shown in \cite{SyYm:3mfd} that in the case of $(2)$ in 
Theorem \ref{thm:graph}, $M$ is homeomorphic to 
an Alexandrov space with nonnegative curvature.

Theorem \ref{thm:graph} determines the possible topological 
type of $M$ if $M$ has not so small diameter. 
In fact, from \cite{CG:collapseI}, every three-dimensional 
graph manifold $M$ has a 
Riemannian metric $g_{\epsilon}$ with sectional curvature 
$|K_{g_{\epsilon}}|\le 1$, $\diam(M,g_{\epsilon})\ge\delta_0$  and 
$\vol(M,g_{\epsilon}) < \epsilon$ for each $\epsilon>0$.

In the bounded curvature case, it follows essentially from 
\cite{CG:collapseII} that  
if a closed three-manifold has a Riemannian metric with $|K|\le 1$ 
whose volume is sufficiently small, then it is a graph
manifold.

The strategy of our proof is as follows: 
We assume $M$ has large diameter which is the essential case.
Applying our previous work \cite{SyYm:3mfd}, we obtain a local 
fiber structure on a neighborhood $B_p$ of each point $p\in M$
over a metric ball $X_p$ in some Alexandrov space 
with curvature bounded below, where $\dim X_p\in\{ 1, 2\}$.
If $\dim X_p=2$, we have a local $S^1$-action on $B_p$. 
If $\dim X_p=1$, we have a (singular) sphere or torus bundle structure
on $B_p$ over the closed interval $X_p$,
and $B_p$ is homeomorphic to one of six compact three-manifolds, which 
will be called {\em cylindrical} if it is homeomorphic to 
either $S^2\times I$ or $T^2\times I$, or {\em cylindrical with a cap} 
if it is homeomorphic to $D^3$, $P^2\tilde\times I$, $S^1\times D^2$ 
or $K^2\tilde\times I$, 
where $\tilde\times$ indicates the twisted product.
Using those local data, we decompose $M$ into two parts as
$M=U_1\cup \hat U_2$, where $U_1$ is a closed domain which 
looks one-dimensional and $\hat U_2$ is one which 
looks two-dimensional, in the Gromov-Hausdorff sense.
More precisely, $U_1$ is defined as the union of all $B_p$
with area of $X_p$ sufficiently small.
Applying the critical point theory for distance functions,
we conclude that each component of $U_1$ is either cylindrical 
or cylindrical with a cap.
We shall construct a local 
$S^1$-action on the remaining piece $\hat U_2$, from which
a graph manifold structure on $M$ is obtained.
To do this,  we need a gluing procedure, which is the main part of the 
present paper. To make the gluing procedure explicit and clear, we give
quantitative descriptions of the local fibering
$B_p\to X_p$ using the geometric properties of fibers in \cite{Ym:collapsing}
and \cite{Ym:conv} over a regular part of $X_p$.
Here the notion of strain radius comes in to control the behavior of the 
regular fibers. This forces us to obtain a sort of compactness of the set 
of regular parts of $X_p$'s  with $B_p$ meeting $\hat U_2$.
This is the reason why in the decomposition $M=U_1\cup \hat U_2$
a neighborhood $B_p$ is included in the 
one-dimensional part $U_1$ even if $\dim X_p=2$ when $X_p$ has a small
area.

The organization of this paper is as follows:
In Section \ref{sec:surf}, we first  establish a uniform lower bound 
on the strain radii of regular parts of $X_p$'s with $B_p$ meeting 
$\hat U_2$, and then provide some basic 
properties of Alexandrov surfaces to show that there are a lot of 
possibilities for the choices of the metric ball $X_p$ 
with boundary having nice geometric properties.
In Section \ref{sec:llocal}, we describe the geometry and topology of 
the local fibering $B_p\to X_p$ in detail.
In Section \ref{sec:decomp}, using these fiber structure 
we have the decomposition $M=U_1\cup \hat U_2$  and 
determine the topology of $U_1$.
In Section \ref{sec:glueI}, we provide a preliminary gluing
argument for the construction of local $S^1$-action on 
$\hat U_2$. The gluing procedure is completed in Section \ref{sec:glueII}.
In Section \ref{sec:thick}, we discuss a thick-thin decomposition 
of a closed orientable Riemannian three-manifold with a lower curvature
bound.

An announcement in Perelman's paper \cite{Pr:surgery}
has recently been come to our attention. He claims that if a 
three-manifold collapses under a local lower sectional curvature bound, 
then it is a graph manifold (Theorem 7.4).  This result also follows from 
the argument in our Theorem \ref{thm:graph} without the extra 
assumption (3) there, since our gluing argument in 
Sections \ref{sec:glueI} and \ref{sec:glueII} is only local 
(see also Section \ref{sec:loc-col}).
The authors do not know his proof of the statement above, up to now.
\par

\begin{ack}
  The authors would like to thank Grisha Perelman
  for correcting our claim in the first draft, 
  on the relation between his Theorem 7.4 and our result,
  as above.
  The second author would like to thank  John Morgan
  and Xiaochun Rong for the discussion at the American Institute of 
  Mathematics,  Palo Alto.
\end{ack}


\section{Strain radii and geometry of Alexandrov surfaces} \label{sec:surf}

We discuss some basic properties of strain radii and metric balls in 
Alexandrov surfaces with curvature bounded below. See \cite{BGP} 
for general facts on Alexandrov spaces.

Let $X$ be an $m$-dimensional complete Alexandrov space 
with curvature bounded below, say, curvature $\ge -1$.
For two points $x, y$ in $X$, a minimal geodesic joining $x$ to $y$
is denoted by $xy$.
The angle between minimal geodesics $xy$ and $xz$ is denoted by $\angle yxz$. 
For a geodesic triangle $\Delta xyz$ in $X$ with vertices $x, y$ and $z$, 
we denote by  $\tilde\angle xyz$ the corresponding angle at $\tilde y$ of
a {\it comparison triangle} 
$\Delta\tilde x\tilde y\tilde z$ for $\Delta xyz$
in the hyperbolic plane of constant curvature $-1$.

For $\delta>0$, 
the {\em $\delta$-regular set} $R_{\delta}(X)$ is defined as the set of 
points $p\in X$ such that there exists 
$m$ pairs of points, $(a_i,b_i)$, $1\le i\le m$, called a 
{\em $\delta$-strainer} at $p$, such that 
\begin{align*}
   \tilde\angle a_ipb_i > \pi - \delta, \quad & 
                \tilde\angle a_ipa_j > \pi/2 - \delta, \\
      \tilde\angle b_ipb_j > \pi/2 - \delta,\quad &
                \tilde\angle a_ipb_j > \pi/2 - \delta,
\end{align*}
for every $i\neq j$. 
The number $\min\,\{ d(a_i,p), d(b_i,p)\,|\, 1\le i\le m\}$ 
is called the {\em length} of the strainer. 
The {\em $\delta$-strain radius} at $p$, denoted by $\delta$-$\strrad(p)$, 
is defined as the supremum of such $r>0$ that there exists a 
$\delta$-strainer at $p$ of length $r$. 
For a closed domain $D$ of $R_{\delta}(X)$, the $\delta$-strain 
radius of $D$, denoted by $\delta$-$\strrad(D)$, is defined as
the infimum of  $\delta$-$\strrad(p)$ when $p$ runs over $D$.
It should be noted that the notion of strain radius is a natural 
generalization of that of injectivity radius for Riemannian manifolds.

For $1\le n\le m$, an $(n,\delta)$-strainer at $p$ is defined by 
$n$ pairs of points, $\{(a_i,b_i)\}$, satisfying the same inequalities 
as above.

For a subset $C$ of $X$,
we denote by $B(C,r)$ or $B(C,r;X)$ the closed metric $r$-ball around 
$C$ and by $S(C,r)$ or $S(C,r;M)$ the metric
$r$-sphere around $C$. For $r<R$, $A(C;r,R)$ denotes the closure 
of $B(C,R)-B(C,r)$.

\begin{lem} \label{lem:str.rad}
For any $m$, $a>0$, $d>0$, $r>0$  and $\delta>0$, there exists a 
positive number
$s=s_m(a,d,r,\delta)$ such that if $B$ is a metric ball in an 
$m$-dimensional complete Alexandrov space $X$ with curvature $\ge -1$
satisfying 
\begin{equation}
    \area(B)\ge a, \quad \diam(B)\le d,
                       \label{eq:geom-bd}
\end{equation}
then the closure $D$ of  $B - B(S_{\delta}(X),r)$
has a definite lower bound for the strain radius$:$
\begin{equation}
  \text{$\delta$-$\strrad(D)\ge s$}. \label{eq:str-rad}
\end{equation}
\end{lem}

\begin{proof}
Certainly we have a positive number $s_X$ depending on $X$ with 
$\delta$-$\strrad(D)\ge s_X$. 
Since the set of all isometry classes of $m$-dimensional 
compact Alexandrov spaces 
satisfying \eqref{eq:geom-bd}
is compact with respect to the Gromov-Hausdorff distance,
this provides a uniform positive lower bound $s_m(a,d,r,\delta)$  for 
all $s_X$.
\end{proof}

Thus the domain $D$ has ``bounded geometry'' in the sense of 
\eqref{eq:str-rad}. This elementary fact is important 
in our gluing argument in Section \ref{sec:glueI}.

The complement $S_{\delta}(X):= X - R_{\delta}(X)$ is called 
the {\em $\delta$-singular set}.
Setting $S_{\delta}(\interior X):=S_{\delta}(X)\cap\interior X$, 
we note that $S_{\delta}(X)=S_{\delta}(\interior X)\cup\partial X$.
Let $ES(\interior X)$ denote the {\em essential singular set} of 
$\interior X$, i.e.,
the set of points $p\in\interior X$ with radius
\[
 {\rm rad}(\Sigma_p):= 
   \min_{\eta\in\Sigma_p} \max_{\xi\in\Sigma_p} d(\xi, \eta) \le \pi/2.
\]

From now on,  we assume $m=2$. Then 
$X$ is known to be a topological two-manifold possibly with boundary.
Moreover $S_{\delta}(\interior X)$ is discrete for any $\delta>0$.

\begin{lem} \label{lem:ess-sing}
For any $p\in X$, $\delta>0$ and $D>0$,
the number of elements of $S_{\delta}(\interior X)\cap B(p,D)$ has a 
uniform upper bound ${\rm Const}(\delta,D)$.
 \par
  In particular we have 
  \[
       \# (ES(\interior X)\cap B(p, D))\le {\rm Const}(D).
  \]
\end{lem}

This follows from an argument similar to Corollary 14.3 of \cite{SyYm:3mfd},
and hence the proof is omitted.

For every $p\in X$, $d_p$ denotes the distance function from $p$.
$d_p$ is called regular at a point $q\neq p$ if there exists 
a $\xi\in\Sigma_q$ such that the directional derivatives of $d_p$
satisfies $d_p'(\xi)>0$.

\begin{lem}[\cite{ST:cut}] \label{lem:cut}
For a fixed $p\in X$, there exists a set $\cal{E}\subset (0,\infty)$ 
of measure zero such that for every $t\in (0,\infty)-\cal E$
\begin{enumerate}
 \item $t$ is a regular value of $d_p ;$ 
 \item $B(p,t)$ is a topological manifold with $($possibly empty$)$
       rectifiable boundary.
\end{enumerate}
\end{lem}

As a consequence of Lemmas \ref{lem:ess-sing} and \ref{lem:cut},
we have 

\begin{cor} \label{cor:collar}
There exists a positive number $\sigma=\sigma(\delta)$  
satisfying the following$:$
For every Alexandrov surface $X$ as above with $\diam(X)>2$, 
for every $p\in X$ and 
for every $t\in [1/2,1]$, there exists $\rho\in (t-10^{-2},t+10^{-2})$
such that 
\begin{enumerate}
 \item $B(p,\rho)$ is a topological manifold$;$
 \item $B(S(p,\rho), \sigma)\cap\interior X \subset R_{\delta}(X);$
 \item $B(S(p,\rho), \sigma)$ is homeomorphic to $S(p,\rho)\times (0,1)$.
\end{enumerate}
\end{cor}

\begin{proof}
The existence of $\sigma$ satisfying (1) and (2) above
is immediate from Lemmas \ref{lem:ess-sing} and \ref{lem:cut}.
For (3), 
it suffices to prove that for every $R>0$ there exists a positive 
constant $\const(R)$ such that if $B(p, R)$ is a topological manifold, 
then the Euler number satisfies  $\chi(B(p,R)) \ge -\const(R)$.
Suppose this does not hold. Then we have a sequence 
$(X_i,p_i)$ of pointed Alexandrov surfaces with curvature 
$\ge -1$ with uniformly bounded $R_i$ such that 
\begin{enumerate}
 \item $B(p_i,R_i)$ is a topological manifold $;$
 \item $\chi(B(p_i,R_i))\to -\infty.$
\end{enumerate}
We may assume that $(X_i,p_i)$ converges to a pointed Alexandrov surface
$(X,p)$. If $\dim X=1$, it is not hard to see that 
$B(p_i,R_i)$ is either a cylinder or a M\"obius band.
If $\dim X=2$, then take an $R$ with $R_i\le R$ for every $i$ and 
choose a regular value $S$ of $d_p$ with
$S>R$  such that $B(p,S)$ is a topological manifold.
Then $B(p_i, S)$ is homeomorphic to $B(p,S)$ by the stability result
(see \cite{Pr:alex2}). This is a contradiction.
\end{proof}

In what follows, we let $\delta^*:=\delta$, which is a sufficiently small 
positive number determined later on in \eqref{eq:delta}.
We also denote the constant $\sigma$ given in Corollary \ref{cor:collar}
by
\begin{equation}
    \sigma^* := \sigma(\delta^*).
                \label{eq:sigma}
\end{equation}



\section{Local structure} \label{sec:llocal}

A local $S^1$-action $\psi$  on a three-manifold $M$ possibly with 
boundary consists of 
an open covering 
$\{ U_{\alpha}\}$ of $M$ and a nontrivial $S^1$-action $\psi_{\alpha}$
on each $U_{\alpha}$ such that
both the actions $\psi_{\alpha}$ and $\psi_{\beta}$
coincide up to orientation on the intersection
$U_{\alpha}\cap U_{\beta}$.
Let $X:=M/S^1$, and $\pi:M\to X$ be the projection. 
$X$ is a topological two-manifold (see \cite{Or:Seif} for instance).

Set
\begin{equation}
  \partial_0 X := \pi(\partial M),\qquad 
   \partial_* X  := \overline{\partial X - \partial_0 X}.
                \label{eq:bdy}
\end{equation}
The fixed point set of $\psi$ coincides with $\partial_* X$.

\begin{lem} \label{lem:action}
If a compact three-manifold $M$ admits a local $S^1$-action with
no singular orbits on $\partial M$, then 
it is a graph manifold.
\end{lem}

\begin{proof}
Note that each component $C$ of $\partial_{*} X$ is a circle.
Take a small collar neighborhood $E(C)$ of $C$ in $X$. Then 
$N(C):=\pi^{-1}(\overline{E(C)})$ is a solid torus.
Setting 
\[
 X_0 := X - \bigcup_C \interior N(C),\qquad   M_0:= \pi^{-1}(X_0),
\]
we have the decomposition
\[
    M = M_0 \cup \Bigl(\bigcup_C N(C)\Bigr),
\]
where $C$ runs over all the components of $\partial_{*} X$.
Since $M_0$ is a Seifert fibered space over $X_0$, 
$M$ is certainly a graph manifold.
\end{proof}

In what follows, let $M$ denote an orientable closed Riemannian 
manifold of dimension three satisfying
\begin{equation}
  K\ge -1,\qquad \vol(M)<\epsilon.
     \label{eq:small-vol}
\end{equation}

We shall determine the geometry and topology of local neighborhoods
of $M$. First we recall the topological structure result for 
such an $M$ when it has uniformly bounded diameter.

\begin{thm}[\cite{SyYm:3mfd}] \label{thm:bdd-diam}
For a given $D>0$, there exists a positive constant 
$\epsilon(D)>0$ satisfying the following$:$
If $M$ satisfies $\diam(M)\le D$ and \eqref{eq:small-vol} with 
$\epsilon\le\epsilon(D)$,
then there exists a $($possibly singular$)$ fibration 
$f:M\to X$, where $X$ is a compact Alexandrov space with curvature 
$\ge -1$ and $\dim X\le 2$. The fiber structure of $M$ can be described
in more detail as follows$:$
\begin{enumerate}
 \item If $\dim X=2$, then $f$ is defined by a local $S^1$-action  on $M$
  with a possible exceptional orbit over a point in $ES(\interior X).$
 \item Let $\dim X=1$. If $X$ is a circle, then $M$ is either a
  sphere-bundle or a torus-bundle over $X$. If $X$ is a closed interval, 
  then $M$ is a gluing of 
  $U$ and $V$ along their boundaries, where $U$ and $V$
  are ones of  $D^3$ and  $P^2\tilde\times I$ or 
  ones of $S^1 \times D^2$ and $K^2\tilde\times I.$
 \item If $\dim X=0$, then a finite cover of $M$ is homeomorphic
   to $S^1\times S^2$, $T^3$, a nilmanifold or
   a simply connected Alexandrov space with nonnegative curvature.
\end{enumerate}
\end{thm}

The case when $X$ is a circle was proved in \cite{Ym:collapsing}, and
the essential part of the case of $\dim X=0$ was proved in  
\cite{FY:fundgp}.

By Lemma \ref{lem:action},
a three-manifold $M$ satisfying one of the conclusions in 
Theorem \ref{thm:bdd-diam}
is a graph manifold except the case when 
$M$ has finite fundamental group and $\dim X=0$.
We also obtain some universal 
positive constants $\delta_0$ and  $\epsilon_0$
such that if  $M$ satisfies $\diam(M) < \delta_0$ and
$\vol(M)<\epsilon_0$, then $M$ is 
homeomorphic to one of the spaces in Theorem \ref{thm:bdd-diam} (3).
Thus Theorem \ref{thm:graph} certainly holds in 
the bounded diameter case. 
Therefore from now we assume that $M$ has large diameter:
\begin{equation}
     \diam(M) \gg  1.
                             \label{eq:diam}
\end{equation}

We now determine the topology of a local 
neighborhood of each point of $M$.

A submersion $f:M\to N$ between Riemannian manifolds is called an 
$\epsilon$-{\it almost Riemannian submersion} if \par
\begin{enumerate}
 \item   the diameter of every fiber of $f$ is less than $\epsilon;$ 
 \item  for every point $p\in M$ and every tangent vector $\xi$ at
        $p$ that is normal to the fiber $f^{-1}(f(p))$, 
   \[
      \left|\frac{|df(\xi)|}{|\xi|} - 1\right| < \epsilon.
   \]
\end{enumerate}
Note that an $\epsilon$-almost Riemannian submersion is a fiber bundle 
map since it is proper.

We denote by  $\tau(\epsilon)$ (resp. $\tau(r|\epsilon)$) a function 
of $\epsilon$ (resp. of $r$ and $\epsilon$) with 
$\lim_{\epsilon\to 0} \tau(\epsilon)=0$
(resp. $\lim_{\epsilon\to 0} \tau(r|\epsilon)=0$ for each fixed $r$).

The following is a quantitative version of Theorem \ref{thm:bdd-diam}
(2) (see \cite{Ym:collapsing}).

\begin{cor} \label{cor:I}
There exists a positive number $\epsilon_1^*$ such that 
if a closed three-manifold $M$ with $K\ge -1$ satisfies 
$d_{GH}(M, I)<\epsilon_1^*$ for some closed interval $I$ of length 
$\ge 1/2$, then there exists a singular fibration 
$f:M\to I$ as in Theorem \ref{thm:bdd-diam} such that 
\begin{enumerate}
 \item the diameter of every fiber of $f$ is less than $\tau(\epsilon_1^*);$
 \item the restriction of $f$ to $I_r$ is a $\tau(r|\epsilon_1^*)$-almost
       Riemannian submersion, where $r>0$ and 
       $I_r:=\{ x\in I\,|\,d(x,\partial I)\ge r\}$.
\end{enumerate}
\end{cor}

Later we shall take $\epsilon_1^*$ such as $\epsilon_1^*\ll \sigma^*$
(see \eqref{eq:stretch}). The final choice of $\epsilon_1^*$ will 
be determined at the end of Section \ref{sec:glueI}.

Let $a^*$ be a positive number such that 
if $B$ is a metric $\rho$-ball with $1/10\le \rho\le 1$ in a complete 
Alexandrov surface $X$ with curvature $\ge -1$ and $\area(B)<a^*$, then 
\begin{equation}
    d_{GH}(B, I) < \epsilon_1^*/2,
                 \label{eq:I}
\end{equation}
for some closed interval $I$.

A surjective map $f:M\to X$ between Alexandrov spaces is called an 
$\epsilon$-{\it almost Lipschitz submersion} if \par
\begin{enumerate}
 \item the diameter of every fiber of $f$ is less than $\epsilon;$ \par
 \item for every $p,q\in M$, if $\theta$ is the infimum of 
        $\angle qpx$ when $x$ runs over $f^{-1}(f(p))$, then 
   \[
      \left|\frac{d(f(p),f(q))}{d(p,q)}-\sin\theta\right| < \epsilon.
   \]
\end{enumerate}
Remark that the notion of $\epsilon$-almost Lipschitz submersion is a 
generalization of 
$\epsilon$-almost Riemannian submersion.
The following result was proved in Theorem 0.2 of \cite{Ym:conv}
(see also Theorem 2.2 of \cite{SyYm:3mfd}).

\begin{thm}[\cite{Ym:conv}] \label{thm:lipsubm}
For given $m$ and $s>0$ there exists $\nu>0$ satisfying the following$:$
Let $X$ be an $m$-dimensional complete Alexandrov space with curvature
$\ge -1$ and with $\delta^*$-$\strrad(X)\ge s$. Then if the Gromov-Hausdorff 
distance between $X$ and a complete Riemannian manifold $M$ with 
$K \ge -1$ is less than $\nu$, then there exists a 
$(\tau(\delta^*)+\tau(s|\nu))$-almost Lipschitz submersion $f:M\to X$
which is a locally trivial bundle map.
\end{thm}

The following is a localized and quantitative version of 
Theorem \ref{thm:bdd-diam}.

\begin{thm} \label{thm:loc-fib}
For every $r>0$, there exists a positive constant 
$\epsilon_0=\epsilon_0(a^*,r,\delta^*)$
satisfying the following$:$
For every  $M$ satisfying \eqref{eq:small-vol} and \eqref{eq:diam} with 
$\epsilon\le\epsilon_0$ and for every $p\in M$,
there exist closed domains $B_p$ and $\hat B_p$ around $p$ 
and 
a pointed complete Alexandrov space $(X,x_0)$ 
with curvature $\ge -1$ and $\dim X\in \{ 1,2\}$
such that 
\begin{enumerate}
 \item $B_p$ and $\hat B_p$ are small perturbations of metric balls 
  around $p$ and, $B(p,1/2)\subset B_p\subset \hat B_p \subset B(p,1)$,
  \quad $\hat B_p - \interior B_p\simeq \partial B_p\times I;$
 \item $B_p$  and $\hat B_p$ have fiber structures over concentric metric 
   balls $X_p\subset \hat X_p$ in $X$ around $x_0;$
 \item $B(x_0,1/2)\subset X_p \subset \hat X_p \subset B(x_0,1)$,\qquad
          $\overline{\hat X_p - X_p}\simeq \partial_0 X_p\times I$,
   where $\partial _0$ denotes the topological boundary.
\end{enumerate}
Moreover the fiber structure on $B_p$ in $(2)$ can be described as follows$:$
Let $\pi_p:(\hat B_p, B_p) \to (\hat X_p, X_p)$ be the fiber projection,
and let
\begin{equation*}
   \partial_* X_p  := \overline{\partial X_p - \partial_0 X_p}.
\end{equation*}
be defined as in \eqref{eq:bdy}.

\begin{proclaim}{\emph{Case} $($A$)$}
     $\dim \hat X_p =1$. 
     $($$\hat X_p$ is a closed interval $I$ in this case$)$.
\end{proclaim}
\begin{enumerate}
 \item[$(a)$] $d_{GH}(\hat B_p, \hat X_p)<\epsilon_1^*$, and the diameter 
              of every fiber of $\pi_p$ is less than 
              $\tau(\epsilon_1^*);$
 \item[$(b)$]  The restriction of $\pi_p$ to $I_r$ is a 
       $\tau(r|\epsilon_1^*)$-almost
       Riemannian submersion$;$
  \item[$(c)$] If $\partial_* X_p$ is empty, then $B_p$ is homeomorphic 
    to either  $I\times S^2$  or $I\times T^2;$
 \item[$(d)$] If $\partial_* X_p$ is nonempty, then $B_p$ is homeomorphic 
   to one of  $D^3$,  $P^2\tilde\times I$ , $S^1\times D^2$ and 
   $K^2\tilde\times I$.
\end{enumerate}
\begin{proclaim}{\emph{Case} $($B$)$}
     $\dim \hat X_p =2$.
\end{proclaim}
\begin{enumerate}
 \item[$(a)$] $d_{GH}(\hat B_p, \hat X_p)<\tau(\epsilon)$, and 
              the length of every fiber of $\pi_p$ is less than 
              $\tau(\epsilon);$
 \item[$(b)$]  $B(B_p,\sigma^*)\subset \hat B_p\subset B(B_p,2\sigma^*);$
 \item[$(c)$]  $B(\partial X_p,2\sigma^*)\cap\interior X\subset
               R_{\delta^*}(X);$
 \item[$(d)$] $B(\partial_0 X_p,2\sigma^*)$ is homeomorphic to 
     $\partial_0 X_p\times (0,1);$
 \item[$(e)$] $D:=\hat X_p - B(S_{\delta^*}(X),r)$ satisfies 
  \begin{enumerate}
   \item[$(i)$] $\delta^*$-$\strrad(D)\ge s$, where 
         $s=s_2(a^*,1,r,\delta^*)$
         is the constant as in Lemma \ref{lem:str.rad}$;$
   \item[$(ii)$] the restriction of $\pi_p$ to $D$ is 
         $(\tau(\delta^*)+\tau(s|\epsilon))$-almost Lipschitz 
         submersion which is an $S^1$-bundle$;$
  \end{enumerate}
  \item[$(f)$] $\pi_p$ gives a local $S^1$-action on $\hat B_p$ whose 
       fixed point set corresponds to $\partial_{*} \hat X_p$,
     where there is a possible exceptional fiber over a point 
     $x\in \hat X_p$ only when $x\in ES(\interior X).$
\end{enumerate}
\end{thm}

\begin{proof}
Suppose the theorem does not hold. Then there exist sequences
$\epsilon_i\to 0$ and 
$M_i$ satisfying \eqref{eq:small-vol}  for $\epsilon_i$  and 
\eqref{eq:diam} such that for some 
$p_i\in M_i$, $B(p_i,1)$ does not contain closed domains 
satisfying the above conclusion. Passing to a subsequence if necessary,
we may assume that $(M_i,p_i)$ converges to a pointed 
complete Alexandrov space $(X,x_0)$ with curvature $\ge -1$.
Observe $1\le \dim X\le 2$. 
In view of Corollary \ref{cor:collar} it is possible to take  
metric balls $Y\subset \hat Y$ of $X$ around $x_0$ which are 
topological manifolds, satisfying
\begin{enumerate}
 \item  $B(x_0,1/2)\subset Y\subset B(Y,\sigma^*)
           \subset\hat Y\subset B(Y,2\sigma^*)\subset B(x_0,1);$
 \item  $B(\partial_0 Y,2\sigma^*)\cap
              \interior X\subset R_{\delta^*}(X);$
 \item  $B(\partial_0 Y,2\sigma^*)$ is contained in a neighborhood of 
        $\partial_0 Y$ homeomorphic to $\partial_0 Y\times (0,1);$
 \item  $\overline{\hat Y - Y} \simeq \partial_0 Y \times I$.
\end{enumerate}
If $\area(\hat Y) < a^*$, then $(\hat Y, Y)$ are Gromov-Hausdorff close to 
some closed intervals $(\hat I, I)$,
 and we put $\hat X_p:= \hat I$, $X_p:= I$ in this case.
If  $\area(\hat Y) \ge a^*$, then we put 
$\hat X_p:=\hat Y$ and $X_p:= Y$.
By Theorem \ref{thm:bdd-diam}, Corollary \ref{cor:I} 
and Theorem \ref{thm:lipsubm} together with
Lemma \ref{lem:str.rad}, we obtain  
closed domains $B_i$ and $\hat B_i$ with
$B(p_i,1/2)\subset B_i \subset \hat B_i \subset B(p_i,1)$ such that 
$B_i$  and $\hat B_i$ are fiber spaces over $Y$  and $\hat Y$ 
respectively as described above
satisfying all the conclusions, which is a contradiction.
\end{proof}


\begin{rem} \label{rem:non-uniq}
\begin{enumerate}
 \item The several geometric properties of  
   $B_p\subset \hat B_p$ and $X_p\subset \hat X_p$ in Theorem
   \ref{thm:loc-fib} will be needed in the gluing argument later on.
 \item We will also need to consider a slight deformation of 
     $X_p\subset \hat X_p$ according to requirements.
\end{enumerate}
\end{rem}

From now on, we put 
\begin{equation}
    r:=\sigma^*/100, \qquad s:=s_2(a^*,1,r,\delta^*).   \label{eq:r}
\end{equation}
Then the constant $\epsilon_0=\epsilon_0(a^*,r,\delta^*)$ in 
Theorem \ref{thm:loc-fib} will become universal
(see the end of Section \ref{sec:glueI}). \par

We now recall 
basic geometric properties of the regular fibers of $\pi_p$. \par

\begin{defn}
For a point $p\in M$ suppose that there is a $(2,\delta^*/2)$-strainer
$\{ (a_j,b_j)\}$ at $p$ of length $\ge s/2$. 
Then the subspace of the tangent space at $p$ generated by 
the directions of minimal geodesics joining $p$ to $a_1$ and $a_2$
is called a {\em horizontal subspace} at $p$.
Let a small circle $F$ in $M$ be given in such a way that 
for every $p\in F$ there is a  $(2,\delta^*/2)$-strainer
at $p$ of length $\ge s/2$.
For $\tau>0$, $F$ is called {\em $\tau$-perpendicular to horizontal subspaces} 
if for each point $p\in F$,
the angle $\theta$ between $F$ and every horizontal subspace at $p$ 
satisfies
\begin{equation*}
               |\theta -\pi/2| < \tau.
\end{equation*}
\end{defn}

\begin{lem}[\cite{Ym:collapsing},\cite{Ym:conv}]  \label{lem:fiber}
Let $\pi_p:\hat B_p\to\hat X_p$  and $D\subset\hat X_p$ 
$($resp. $I_r\subset\hat X_p$$)$ be as in Theorem \ref{thm:loc-fib}.
For every $x\in D$ $($resp. $x\in I_r$$)$ and $q\in \pi_p^{-1}(x)$, 
the following holds:
\begin{enumerate}
 \item For every  $q'$with $d(q,q')\ge r$, the angle $\theta$
       between $\pi_p^{-1}(x)$ and every minimal geodesic joining
       $q$ to $q'$ satisfies
   \begin{equation*}
          |\theta -\pi/2| < \tau(\delta^*)+\tau(s|\epsilon)\quad 
                (resp.\,\, \tau(r|\epsilon_1^*));
   \end{equation*}
 \item The fiber $\pi_p^{-1}(x)$ is 
       $(\tau(\delta^*)+\tau(s|\epsilon))$-perpendicular to horizontal 
       subspaces.
\end{enumerate}
\end{lem}


\section{Decomposition} \label{sec:decomp}

Let $M$ satisfy $\eqref{eq:small-vol}$ with $\epsilon\le\epsilon_0$.
Take points $p_1, p_2, \ldots,$ of $M$
such that the collections $\{ B_{p_i}\}$ and  $\{\hat B_{p_i}\}$ 
given by Theorem \ref{thm:loc-fib} are finite coverings of $M$.
We may assume that 
 \begin{align}
   & d(p_i, p_j)\ge 1/10  \quad\text{for every $i\neq j;$}  
                   \label{eq:apart}  \\
   &\text{$\{ B_{p_i}(1/10) \}$ covers $M$,} \label{eq:cover}
 \end{align}
where  
$B_{p_i}(1/10):= \{ x\in B_{p_i}\,|\,d(x,\partial B_{p_i})\ge 1/10\}$.
By the Bishop-Gromov volume comparison theorem, 
we may assume that the maximal number of $\hat B_{p_i}$'s having
nonempty intersection is uniformly bounded above by
a universal constant $Q$  not depending on $M$.
Let $X_{p_i}\subset\hat X_{p_i}$ be chosen as in Theorem \ref{thm:loc-fib}
for $B_{p_i}\subset\hat B_{p_i}$.

For simplicity, we put 
 \[ B_i:=B_{p_i},\quad \hat B_i:= \hat B_{p_i},\quad
   X_i:=X_{p_i}, \quad  \hat X_i:= \hat X_{p_i}.
 \]
If $\dim X_i=2$, then there exists a local 
$S^1$-action $\psi_i$ on $B_i$ such that 
$B_i/\psi_i \simeq X_i$, where 
$\pi_i:=\pi_{p_i}:\hat B_{i}\to \hat X_{i}$ is the projection.
\par

For each $j\in \{ 1,2\}$, 
let $\cal I_j$ denote the set of all $i$ with $\dim X_i=j$,
and consider 
\[
    U_j := \bigcup_{i\in \cal I_j} B_i.
\]

Let $\cal B_j:=\{ B_i\,|\,\dim X_i=j\}$, $j\in\{ 1,2\}$.
By Theorem \ref{thm:loc-fib}, each element of $\cal B_1$ is either 
cylindrical or cylindrical with a cap (see Introduction).

\begin{lem}\label{lem:dim1}
Each component of $U_1$ is homeomorphic to one of 
$D^3$, $P^2\tilde\times I$, $S^2\times I$, $S^1\times D^2$, 
$K^2\tilde\times I$ and $T^2\times I$  unless $U_1=M$.

If $U_1=M$, then $M$ is  homeomorphic to one of the spaces in 
Theorem \ref{thm:bdd-diam} (2).
\end{lem}

\begin{proof}
Slightly enlarging closed domains $B_i$ in $\mathcal B_1$ if necessary, 
we may assume 
that any two $B_i, B_j$ in $\cal B_1$ has intersection 
$B_i\cap B_j$ which is either empty or else having diameter $>\sigma^*$.
Since  $U_1$ is a part of $M$ which looks one-dimensional
in the Gromov-Hausdorff sense, it follows from \eqref{eq:apart} 
that $U_1$ is a manifold. 
Suppose that $B_{i}\cap B_{j}$ is nonempty for two domains 
$B_i, B_j$ in $\cal B_1$.
In the argument below, we may assume that 
$B_{i}\not\subset B_{j}$ and  $B_{j}\not\subset B_{i}$.
Since $B_i$ and $B_j$ are either cylindrical or cylindrical with a cap,
it follows from \eqref{eq:diam} that at least one of $B_i$ and $B_j$, 
say $B_j$, has disconnected boundary.
Consider the distance function $d_{p_i}$, where
$p_i$ is the reference point of $B_i$. 
Letting $F$ denote $S^2$ or $T^2$, we know that $B_i$ and $B_j$ have 
$F$-fiber structures over $I$, which is singular at the top of the cap. 
By Lemma \ref{lem:fiber}, one can construct a gradient-like vector field
$V_i$ for $d_{p_i}$ on a neighborhood of $\overline{B_j-B_i}$
whose flow curves are transversal to every fiber of $B_j$ lying on a
neighborhood of $\overline{B_j-B_i}$.
In view of \eqref{eq:diam}, it follows that $B_i\cup B_j$ is 
homomorphic to $B_i$. Repeating the argument finitely 
many times, we obtain the conclusion of the lemma.
\end{proof}

From now on, we assume $U_1\neq M$, and 
consider the decomposition of $M$ 
\[
       M = U_1\cup \hat U_2,
\]
where $\hat U_2$ denotes the closure of $M-U_1$.

For every fixed component $L$ of $\partial U_1$, there exists a 
unique $B_{\ell}\in\cal B_1$ such that a component of $\partial B_{\ell}$ 
coincides with $L$. Let $B^L_{1},\ldots,B^L_{n}$ denote the set of 
all elements of $\cal B_2$ such that $B^L_{i}(1/10)$ 
meets $L$, $1\le i\le n$. 
Since $\diam(L)<\tau(\epsilon_1^*)$, 
$B^L_{i}$ contains $L$. 
Let $L_i$ be the unique component of $\partial B^L_i$ meeting 
$B_{\ell}$.
Since the domain bounded by $L$ and $L_i$ is one-dimensional
in the Gromov-Hausdorff sense, it follows from the curvature condition
that  $B^L_{1},\ldots,B^L_{n}$ lie in a linear order and 
for every $1\le i \neq j\le n$
\begin{equation}
     d(\partial B^L_i, \partial B^L_j) \ge 1/20,
                 \label{eq:biggest}
\end{equation}
and we may assume that  
\begin{equation}
     B^L_n \supset L_i.
                 \label{eq:biggest}
\end{equation}
We denote by  $\tilde U_2$ the union of all $B_j$ in $\cal B_2$
which does not intersect $U_1$. Obviously we have 
\[
   U_2 =\tilde U_2\cup\Bigl(\bigcup_L (B^L_1\cup\cdots\cup B^L_n)\Bigr),
\]
where $L$ runs over all the components of $\partial U_1$.

\begin{lem}\label{lem:U_2}
        $U_2$ is homeomorphic to $\hat U_2$.
\end{lem}

\begin{proof}
Fix a component $L$ of $\partial U_1$ again.
Let $\hat B^L_n\supset B^L_n$,  $\hat X^L_n\supset X^L_n$
and $\pi^L_n:(\hat B^L_n,B^L_n) \to (\hat X^L_n, X^L_n)$ be the 
orbit projection as in Theorem \ref{thm:loc-fib}.
Let $\hat L_n$ denote the component of $\partial \hat B^L_n$ 
corresponding to $L_n$, and let
$U$ be the domain bounded by $L$ and $L_n$.
Take a point $x\in\hat U_2$ with $d(x, L)\ge 1$ and a point
$y\in\hat L_n$.
Since every fiber of  $\pi_n^L$ meeting $U$ has diameter $<\tau(\epsilon_1^*)$,
it follows that for every $z\in U$
\begin{equation}
     \tilde\angle x z y > \pi - \tau(\sigma^*|\epsilon_1^*).
                    \label{eq:stretch}
\end{equation}
Let $V$ be a gradient-like vector field for $d_{x}$ defined on a neighborhood
of $U$.

\begin{ass} \label{ass:transv}
The flow curves of $V$ are transversal to both $L$ and $L_n$.
\end{ass}

\begin{proof}
Since the transversality to $L$ is immediate from Lemma \ref{lem:fiber} (1),
it suffices to check the transversality to $L_n$.
For every $p\in L_n$, let $\phi_p(t)$ be the flow curve of $V$ 
with $\phi_p(0)=p$. Put $q:=\exp_p \sigma^* V(p)$, $\bar p:=\pi^L_n(p)$
and denote by $\bar\xi$ the direction at $\bar p$ defined by a
minimal geodesic to $\pi^L_n(q)$.
In a way similar to Lemma 4.6 of \cite{Ym:conv},
we have 
\[
 d(\pi^L_n(\phi_p(t)),\exp_{\bar p} t\bar\xi)
     < t (\tau(\delta^*) + \tau(s|\epsilon)),
\]
for every sufficiently small $t>0$.
This implies that $\phi_p(t)$ makes an angle with $L_n$ 
uniformly bounded away from zero.
\end{proof}

Now in view of \eqref{eq:biggest}, Assertion \ref{ass:transv} implies 
$U_2\simeq \hat U_2$.
\end{proof}

The proof of the following lemma is deferred to Sections \ref{sec:glueI} 
and \ref{sec:glueII}.

\begin{lem}\label{lem:glue}
There exists a local $S^1$-action defined on 
$U_2$, and hence on $\hat U_2$.
\end{lem}

\begin{proof}[Proof of Theorem \ref{thm:graph} assuming Lemma
\ref{lem:glue}]
Note that each component of $\partial \hat U_2$ is homeomorphic to 
$S^2$ or $T^2$.
For each component $L$ of $\partial\hat U_2$, 
let $W(L)$ be the component of $U_1$ containing $L$.
Suppose first that $L$ is homeomorphic to $S^2$.
Then one of the following holds:
\begin{enumerate}
 \item $\partial W(L)$ is connected and $W(L)$ is homeomorphic to either 
   $D^3$ or $P^2\tilde\times I;$
 \item $W(L)$ is homeomorphic to $S^2\times I$, and the other component 
   of $\partial W(L)$ is another component of $\hat U_2$.
\end{enumerate}
Now consider the union
\[
  V:= \hat U_2\cup\Bigl(\bigcup_{L} W(L)\Bigr),
 \]
where $L$ runs over all the components of $\hat U_2$ 
homeomorphic to $S^2$. 
Since an $S^1$-action on $S^2$ is essentially by rotation, 
the local $S^1$ action on $\hat U_2$ extends to a 
local $S^1$-action on $V$ such that 
the orbit space $W(L)/S^1$ is a disk whose singular locus is
one of 
\begin{enumerate}
 \item an interval on the boundary of $W(L)/S^1$ ( the case of 
       $W(L)\simeq D^3$);
 \item the union of an interval on the boundary of $W(L)/S^1$ and 
       a  point in $\interior W(L)/S^1$ of type $(2,1)$-singularity
      (the case of $W(L)\simeq P^2\tilde\times I$);
 \item the disjoint union of two intervals on the boundary of $W(L)/S^1$ 
      (the case of $W(L)\simeq S^2\times I$).
\end{enumerate}
Note also that each component of $\partial V$ is homeomorphic to $T^2$
and having no singular orbits. Therefore Lemma \ref{lem:action}
implies that $V$ is a graph manifold.
From construction, for each component $L$ of $\partial V$, 
one of the following holds:
\begin{enumerate}
 \item[(a)] $\partial W(L)$ is connected and $W(L)$ is homeomorphic to 
   either $S^1\times D^2$ 
   or  $K^2\tilde\times I;$
 \item[(b)]$W(L)$ is homeomorphic to $T^2\times I$, and the other component of 
   $\partial W(L)$ is another component of $\partial V$.
\end{enumerate}
Thus  $M$ is a graph manifold. 
\end{proof}
\par\medskip


\section{Gluing} \label{sec:glueI}

Let $B\subset \hat B$ be closed domains in a closed orientable
three-manifold $M$
with sectional curvature $K\ge -1$, and let 
$X\subset \hat X$ be concentric closed metric balls of radii $t<\hat t$ 
in a two-dimensional complete Alexandrov space $Z$ with 
curvature $\ge -1$. 
Assume that
\begin{enumerate}
 \item the Gromov-Hausdorff distance between $\hat B$ and $\hat X$
    (resp. $B$ and $X$) is sufficiently small;
 \item $X\subset \hat X$ satisfy the conclusion of Case (B) in 
    Theorem \ref{thm:loc-fib};
 \item $\area(X)\ge a^*;$
 \item $1/10\le t\le t + s\le \hat t\le 1$,
\end{enumerate}
where $s$ is as in \eqref{eq:r}.
Let $D:=X - B(S_{\delta^*}(Z), r)$, $\hat D:=\hat X - B(S_{\delta^*}(Z), r)$.
Note that $\delta^*$-$\strrad(\hat D)\ge s$.
Applying Theorem \ref{thm:lipsubm}, we have closed domains $\hat N$ and $N$ 
of $\hat B$ and $B$ respectively, and an almost Lipschitz submersion
$\pi:(\hat N, N)\to (\hat D, D)$, which is an $S^1$-bundle.

First we need to establish the uniform boundedness of length ratio for 
the fibers of $\pi:N\to D$.

\begin{lem} \label{lem:length}
There exists a $\zeta=\zeta(a^*,s^*,\delta^*)>0$ 
such that the following holds$:$
Suppose that 
\begin{enumerate}
 \item for every $x\in \hat D$ the length $\ell(x)$ of the fiber  
       $\pi^{-1}(x)$ is less than $\zeta;$
 \item for every $x\in D$and  $p\in \pi^{-1}(x)$, 
       letting $\theta(p)$ denote 
       the angle between $\pi^{-1}(x)$ and 
       a horizontal subspace at $p$, we have 
    \begin{equation*}
       |\theta(p)-\pi/2| < \zeta.
     \end{equation*}
\end{enumerate}
Then 
 \begin{equation*}
        c^{-1} < \frac{\ell(x)}{\ell(y)} < c,
 \end{equation*}
for every $x,y\in D$, where $c=c(a^*,s)$ is a uniform 
positive constant.
\end{lem}

\begin{proof}
Suppose the lemma does not hold. Then we have a sequence 
$\pi_i:(\hat N_i, N_i)\to (\hat D_i, D_i)$ of $S^1$-bundles satisfying 
the assumptions of the lemma for $\zeta_i$ with $\lim \zeta_i=0$
such that 
$\frac{\ell(x_i)}{\ell(y_i)}\to \infty$, where 
$\pi_i^{-1}(x_i)$ (resp. $\pi_i^{-1}(x_i)$)  has the maximal 
(resp. the minimal) length among all the fibers of $\pi_i$ over $D_i$. 
Note that 
\begin{equation*}
   N_i\simeq  
    \begin{cases}
      D_i\times S^1 \quad & \text{if $D_i$ is orientable}, \\
      D_i\tilde\times S^1 \quad & \text{if $D_i$ is non-orientable}.
     \end{cases}
\end{equation*}
Take a finite covering $(\hat E_i, E_i) \to (\hat N_i, N_i)$ 
along fibers such that 
the length $\hat \ell(x_i)$ of $\hat\pi_i^{-1}(x_i)$ satisfies 
$1< \hat \ell(x_i) < 2$, where $\hat \pi_i:\hat E_i \to \hat D_i$ 
is the natural projection. We may assume that 
$\hat \pi_i:(\hat E_i, E_i) \to (\hat D_i, D_i)$ 
converges to a Lipschitz map  $\hat\pi:(\hat E, E) \to (\hat D, D)$. 
Note that $D\subset \hat D$ are closed domain in some complete 
Alexandrov surface with curvature $\ge -1$ satisfying
$\delta^*$-$\strrad(\hat D)\ge s$.
Let $x$ be the limit of $x_i$, and  $F:=\hat\pi^{-1}(x)$.
Choose a $\delta^*$-strainer $\{ (a_j,b_j)\}_{j=1,2}$ at $x$ of length $s$.
For every $p\in F$, one can take
$q_j\in\pi^{-1}(a_j)$ and $r_j\in\pi^{-1}(b_j)$
such that  $\{ (q_j,r_j)\}$ is a strainer at $p$
and that $pq_j$ and $pr_j$ are almost perpendicular to $F$. 
Since $F$ has a positive 
diameter, this implies $\dim \hat E = 3$.
Thus $\hat E_i$ does not collapse. Since $E\subset R_{\delta^*}(\hat E)$, 
it follows from Theorem \ref{thm:lipsubm} that $E_i$ is almost isometric 
to $E$ for large $i$ in the sense that there is a bi-Lipschitz
homeomorphism $f_i:E_i \to E$ such that the Lipschitz constants
of $f_i$ and $f_i^{-1}$ are close to one.
Note that the length of shortest nonzero homotopic loops in $E$
has a definite positive lower bound, while 
the length of $\hat\pi_i^{-1}(y_i)$ converges to zero.
This is a contradiction. 
\end{proof}

Let $\pi:(\hat B, B)\to (\hat X, X)$, $(\hat D, D)\subset (\hat B,B)$,
$\pi:(\hat N, N)\to (\hat D, D)$ and $\hat B\subset Z$  be as in 
Lemma \ref{lem:length}.
Let $N_1$ be the set of points  $p$ of $B$
such that there is a $(2, 2\delta^*)$-strainer at $p$ of length $\ge s/2$.
Note $N\subset N_1$.

\begin{lem} \label{lem:N_1}
Every fiber $F$ of $\pi$ contained in $N_1$ is 
$(\tau(2\delta^*)+\tau(s/2|\epsilon))$-perpendicular to
horizontal subspaces.
\end{lem}

\begin{proof}
It is obvious that $\pi(N_1)$ is contained in 
$R_{3\delta^*}(Z)$ and having $3\delta^*$-strainer of length $\ge s/3$.
Thus  Theorem \ref{thm:lipsubm} together
with Lemma \ref{lem:fiber} (2) yields the conclusion.
\end{proof}

Next we consider a gluing situation.
Let $B$ and  $B'\subset \hat B'$ be closed domains in  $M$,
and let 
$\pi:B\to X$ and $\pi':(\hat B', B') \to (\hat X', X')$ 
be the orbit maps of local $S^1$-actions,
where we assume that $X$ is only a topological two-manifold.
Here we consider $\pi:B\to X$ as a result of gluing of several 
local $S^1$-actions 
$\{ \pi_i:B_i\to X_i\}$. Note those $S^1$-actions are neither 
isometric nor our gluing will be through isometric actions.
This is the reason why we assume $X$ to be only a topological 
two-manifold.
On the other hand, we assume 
$(\hat X', X')$,  $(\hat D', D')$ and 
$\pi':(\hat N', N')\to (\hat D', D')$ 
are as in Lemma \ref{lem:length},
where  
$D':=X' - B(S_{\delta^*}(Z'), r)$, 
$\hat D':=\hat X' - B(S_{\delta^*}(Z'), r)$, 
and $Z'$ is a two-dimensional complete Alexandrov space with 
curvature $\ge -1$ containing $(\hat X', X')$.

Let $N_1\subset B$, $N_1'\subset B'$ be defined as in 
Lemma \ref{lem:N_1}, and suppose $B\cap B'$ is nonempty.

In the sequel, for a closed domain $A$ of $M$ we denote
by $\underline{A}$ a small perturbation of $A$. 

\begin{lem} \label{lem:glue1}
For a given positive number $\nu$ there exist 
$\delta^*=\delta^*(\nu)>0$ and $\zeta=\zeta(\nu)>0$ such that 
the following holds$:$
Let $\pi:B\to X$,  $\pi':(\hat B', B') \to (\hat X', X')$ 
be as above satisfying
\begin{enumerate}
 \item[$($a$)$] every fiber $F$ of $\pi$ contained in $N_1$ and $F'$ 
   of $\pi'$ contained in $N'$  are $\zeta$-perpendicular to horizontal 
   subspaces$;$
 \item[$($b$)$] any two orbits of $\pi$ in $N$ with distance 
   $\le 1$ have length ratio uniformly bounded as in Lemma \ref{lem:length}.
\end{enumerate}
We also assume that
\begin{equation}
     \pi'(\partial B'\cap B)\subset D'.
                      \label{eq:bdy-reg}
\end{equation}
Then we have a local $S^1$-action $\psi''$ on a small perturbation 
$\underline{B\cup B'}$ of $B\cup B'$ 
and a topological two-manifold $X''$
with the orbit projection $\pi'':\underline{B\cup B'}\to X''$
satisfying
\begin{enumerate}
\item $\underline{B\cup B'}$ is a manifold with boundary, and 
 \begin{equation*}   
   (B\cup B')(10\ell')\subset \underline{B\cup B'}\subset B(B\cup B',10\ell'),
 \end{equation*}
   where $\ell'$ denotes the maximal length of fibers of $\pi'$
   meeting $\partial B'\cap B$, 
   and $(B\cup B')(10\ell')=\{ x\in B\cup B'\,|\,d(x,\partial(B\cup B'))\ge
       10\ell'\};$
 \item 
  \begin{equation*}
    \psi''= \begin{cases}
      \psi &  \quad \text{on \,
              $\underline{B\cup B'}-B(B',10\ell')$}  \\
      \psi' & \quad \text{on \,
              $B'$};
        \end{cases}
  \end{equation*}
 \item each orbit of $\psi''$ has length  $<2\ell''$, where
    $\ell''$ denotes the maximal length of all fibers of $\pi$ and $\pi'$
    intersecting $10\ell'$-neighborhood of $\partial B'\cap B;$
 \item every fiber of $\pi''$ contained in $N_1''$ is 
       $\nu$-perpendicular to horizontal subspaces, where $N_1''$ is 
       the set of points  $p$ of $\underline{B\cup B'}$ such that 
       there is a $(2, 2\delta^*)$-strainer at $p$ of length $\ge s/2$.
\end{enumerate}
\end{lem}

\begin{rem} 
Under the situation of Lemma \ref{lem:glue1}, 
if both $(\hat B, B)$ and $(\hat B', B')$ are as in Lemma \ref{lem:length},
then by Lemmas \ref{lem:length} and \ref{lem:N_1}, $B\to X$ and 
$\pi':(\hat B', B')\to (\hat X', X')$ satisfy the assumptions of 
Lemma \ref{lem:glue1} except \eqref{eq:bdy-reg} if 
$\tau(2\delta^*)+\tau(s/2|\epsilon) < \zeta$ which is 
realized by $\delta^*\ll 1$ and $\epsilon\ll \delta^*$.
Note that the proof of  Lemma \ref{lem:length} goes through for $N_1$ 
as well in place of $N$.
\end{rem}

For the proof of Lemma \ref{lem:glue1}, we need a sublemma.

Let $A$ be a small neighborhood of $\pi'(\partial B'\cap B)$
in $\pi'(\partial B')$, and let $C$ be the closure of the intersection 
of $\interior X$ with the boundary of the $10\ell'$-neighborhood of 
$\pi(B\cap B')$.
Slightly perturbing $A$ and $C$ if necessary, we may assume that 
both are one-manifolds. 

Fix any $x, \hat x\in A$ with 
$10\ell'\le d(x,\hat x) \le 20\ell'$.
Taking a nearest point $z$ of $\pi'(\pi^{-1}(C))$ from $x$,
choose any point $y\in \pi((\pi')^{-1}(z))$.
Similarly we choose $\hat y\in C$ for $\hat x$.
Put $F':=(\pi')^{-1}(x)$, $\hat F':=(\pi')^{-1}(\hat x)$,
$F:=(\pi)^{-1}(y)$ and $\hat F:=(\pi)^{-1}(\hat y)$. 

\begin{slem} \label{slem:isotopy}
Under the situation above, there exist $\delta^*=\delta^*(\nu)>0$ 
and $\zeta=\zeta(\nu)>0$ satisfying the following$:$
\begin{enumerate}
 \item \begin{equation*}
          \left|\frac{\ell(F)}{\ell(F')} -1\right|
                      < \nu,
       \end{equation*}
    where $\ell(F)$ denotes the length of $F$.
 \item There exists an  annulus 
       $E$  $($resp. $\hat E$$)$ in $M$
       equipped with an $S^1$-fiber structure via a bi-Lipschitz homeomorphism 
       $h:[0,1]\times S^1\to E$ $($resp. $\hat h:[0,1]\times S^1\to\hat E$$)$
      such that 
       \begin{enumerate}
         \item $F=h(0\times S^1)$ and $F'=h(1\times S^1)$
           $($resp.  $\hat F=\hat h(0\times S^1)$ and 
             $\hat F'= \hat h(1\times S^1)$$);$
         \item for each $t\in [0,1]$, $h(t\times S^1)$
           $($resp. $\hat h(t\times S^1)$ $)$
           is $\nu$-perpendicular to horizontal subspaces.
       \end{enumerate}
 \item Let $[x, \hat x]$ and $[y, \hat y]$ be the subarcs
       of $A$ and $C$ respectively, and let $T$ be the union of 
       $(\pi')^{-1}([x, \hat x])$, $\pi^{-1}([y, \hat y])$, 
       $E$ and $\hat E$.
      Then the domain $D$ bounded by $T$ has 
      an $S^1$-fiber structure via a bi-Lipschitz homeomorphism 
      $k:D^2\times S^1\to D$ such that 
       \begin{enumerate}
         \item for each $x\in\partial D^2$, $k(x\times S^1)$
            coincides with a fiber on $T;$
         \item for each $x\in D^2$, $k(x\times S^1)$ 
            is $\nu$-perpendicular to horizontal subspaces.
       \end{enumerate}
\end{enumerate}
\end{slem}

\begin{figure}[htbp]
  \centering
  \includegraphics[scale=0.8]{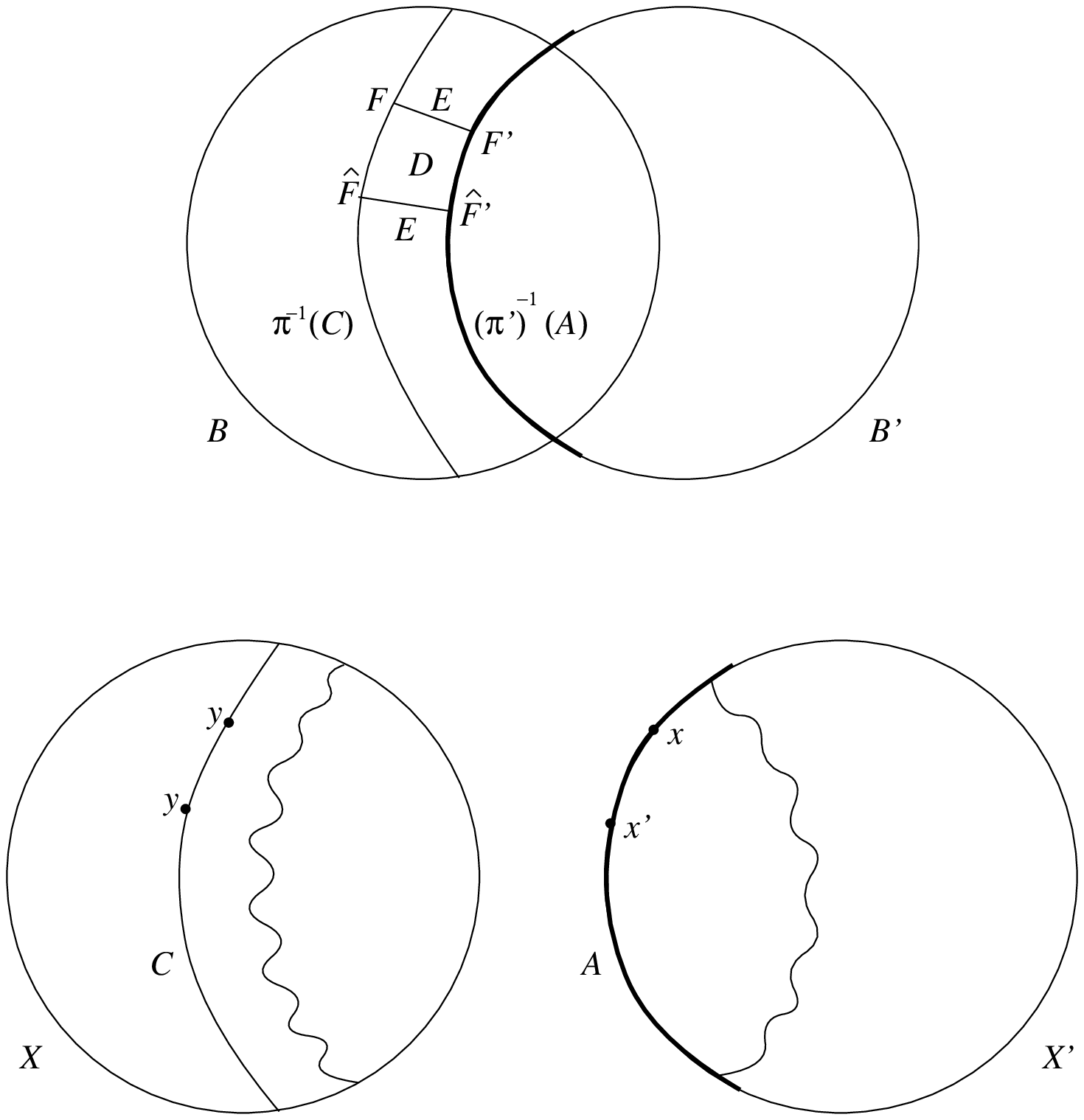}
  \caption{}
  \label{fig:S}
\end{figure}

\begin{proof}
We prove it by contradiction. 
If the conclusion does not hold, we would have a sequence of 
closed three-manifolds $M_i$ with $K\ge -1$ for which 
there are $\pi_i:B_i\to X_i$, $\pi_i':(\hat B_i', B_i')\to (\hat X_i', X_i')$
satisfying the assumptions of Sublemma \ref{slem:isotopy}
for $\delta_i^*\to 0$ and $\zeta_i\to 0$, but not satisfying the conclusions 
for $\zeta$.
Let $A_i\subset \partial X_i'$,
$C_i\subset X_i$, $x_i, \hat x_i\in A_i$ and
$y_i, \hat y_i\in C_i$  be defined as above. 
In particular, $10\ell_i'\le d(x_i, \hat x_i)\le 20\ell_i'$, 
where $\ell_i'$ is defined in a way similar to $\ell'$
(see Lemma \ref{lem:glue1} (1)).
Put  $F_i':=(\pi_i')^{-1}(x_i)$,
$\hat F_i':=(\pi_i')^{-1}(\hat x_i)$,
$F_i:=(\pi_i)^{-1}(y_i)$ and  $\hat F_i:=(\pi_i)^{-1}(\hat y_i)$.
Let $\ell_i(y_i)$ and $\ell_i'(x_i)$ denote
the length of $F_i$ and $F_i'$ respectively.
Passing to a subsequence, we may assume that 
$(\frac{1}{\ell_i'(x_i)} M_i,x_i)$ converges to a pointed space
$(W,w_0)$, where $W$ is a complete Alexandrov space with nonnegative
curvature. From assumption, we see that 
$W$ is actually isometric to $\R^2\times S^1_1$, where $S^1_1$ denotes 
the circle of length $1$.
Thus for any fixed $R\gg 1$, $B(x_i, R;\frac{1}{\ell_i'(x_i)} M_i)$
is almost isometric to $B(w_0, R)$.
This together with the condition (a) of Lemma \ref{lem:glue1}
 implies that $\ell_i'(x_i)$ and the length of any orbit of 
$\pi_i$ nearby $(\pi_i')^{-1}(x_i)$ are comparable 
in the sense of Lemma \ref{lem:length}.
Then by \eqref{eq:bdy-reg}, $\ell_i'(x_i)$ and 
$\ell_i(y_i)$ are comparable. The above convergence then 
yields $\frac{\ell_i(y_i)}{\ell_i'(x_i)}\to 1$, which proves (1).

Let $F$, $\hat F$,   $F'$ and  $\hat F'$ be the limits of 
$F_i$, $\hat F_i$,   $F_i'$ and  $\hat F_i'$ respectively
under the above convergence.
Since $F$ and $F'$ (resp. $\hat F$ and $\hat F'$) can be joined by 
one-parameter family of parallel
circles, say $E$ (resp. say $\hat E$) of length $1$, 
$F_i$ and $F_i'$ can be joined by one-parameter family, say $E_i$
(resp. say $\hat E_i$)  of 
circles each of which is $\nu$-perpendicular to horizontal subspaces for 
sufficiently large $i$.
Let 
$\varphi_i: B(w_0, R)\to  B(x_i, R;\frac{1}{\ell_i'(x_i)} M_i)$
be an almost isometry.
Note that the closed domain $D_i$ bounded by 
$\pi_i^{-1}([y_i,\hat y_i])$, $(\pi_i')^{-1}([x_i,\hat x_i])$,
$E_i$ and $\hat E_i$  is mapped via $\varphi_i^{-1}$ onto 
a domain $D$ bounded by $E$, $\hat E$, $\pi^{-1}([x,\hat x])$
and $\pi^{-1}([y,\hat y])$. 
Note that $\varphi_i$ maps horizontal subspaces to horizontal subspaces
(see \cite{Ym:collapsing} for the details).
Since $D$ is isometric to a product $H\times S^1_1$ for a rectangle $H$,
this gives a compatible $S^1$-fiber 
structure on $D_i$ each of whose fibers is $\nu$-perpendicular 
to horizontal subspaces. This is a contradiction.
\end{proof}

\begin{proof}[Proof of Lemma \ref{lem:glue1}]
We shall carry out the required gluing procedure on each component, say $U$,
of $B\cap B'$.
Let $A_0$ be any component of $A\cap\pi'(U)$, and take consecutive points
$x_1,\ldots, x_N$ of $A_0$
with $10\ell'\le d(x_{\alpha}, x_{\alpha+1})\le 20\ell'$ for 
each $1\le \alpha \le N-1$.

First consider

\begin{proclaim}{\emph{Case} (A)}
    $\partial B$ does not meet $\partial B'$ on $U$.
\end{proclaim}

\begin{figure}[htbp]
  \centering
  \includegraphics[scale=0.8]{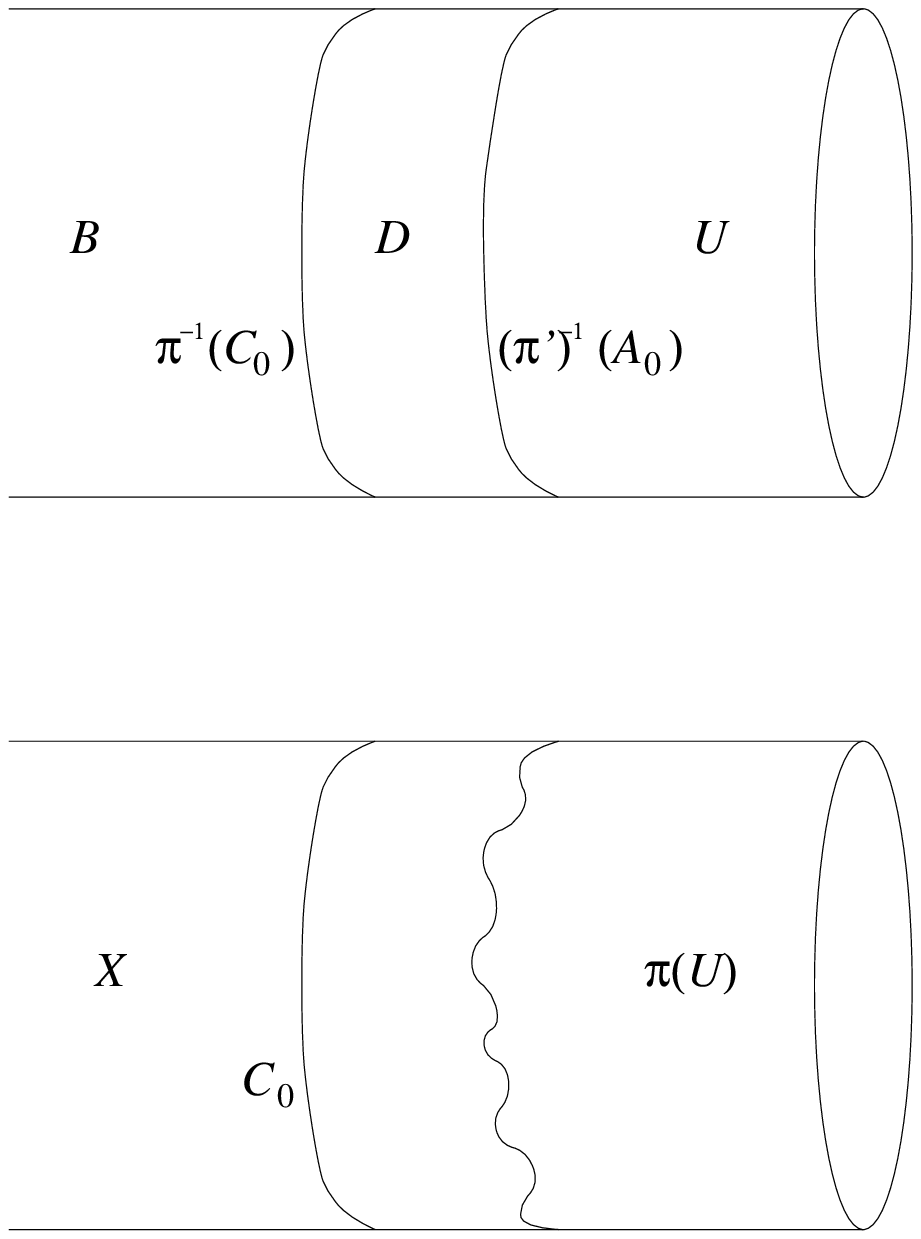}
  \caption{Case (A)}
  \label{fig:A}
\end{figure}
%
In this case, both $A_0$ and the component $C_0$ of $C$ 
corresponding to $A_0$
are circles. 
Applying Sublemma \ref{slem:isotopy} to $x:=x_{\alpha}$ and 
$\hat x:=x_{\alpha+1}$, 
we obtain a closed domain $D$ bounded by $\pi^{-1}(C_0)$ and $(\pi')^{-1}(A_0)$
having an $S^1$-bundle structure via a bi-Lipschitz homeomorphism 
$k:(I\times S^1)\times S^1\to D$ such that 
\begin{enumerate}
 \item for each $x\in \partial I\times S^1$, $k(x\times S^1)$
      coincides with a fiber on 
      $\pi^{-1}(C_0)\cup (\pi')^{-1}(A_0);$
 \item for each $x\in I\times S^1$, $k(x\times S^1)$ is 
      $\nu$-perpendicular to horizontal subspaces.
\end{enumerate}

\begin{proclaim}{\emph{Case} (B)}
  $\partial B$ meets $\partial B'$ on $U$.
\end{proclaim}

\begin{figure}[htbp]
  \centering
  \includegraphics[scale=0.8]{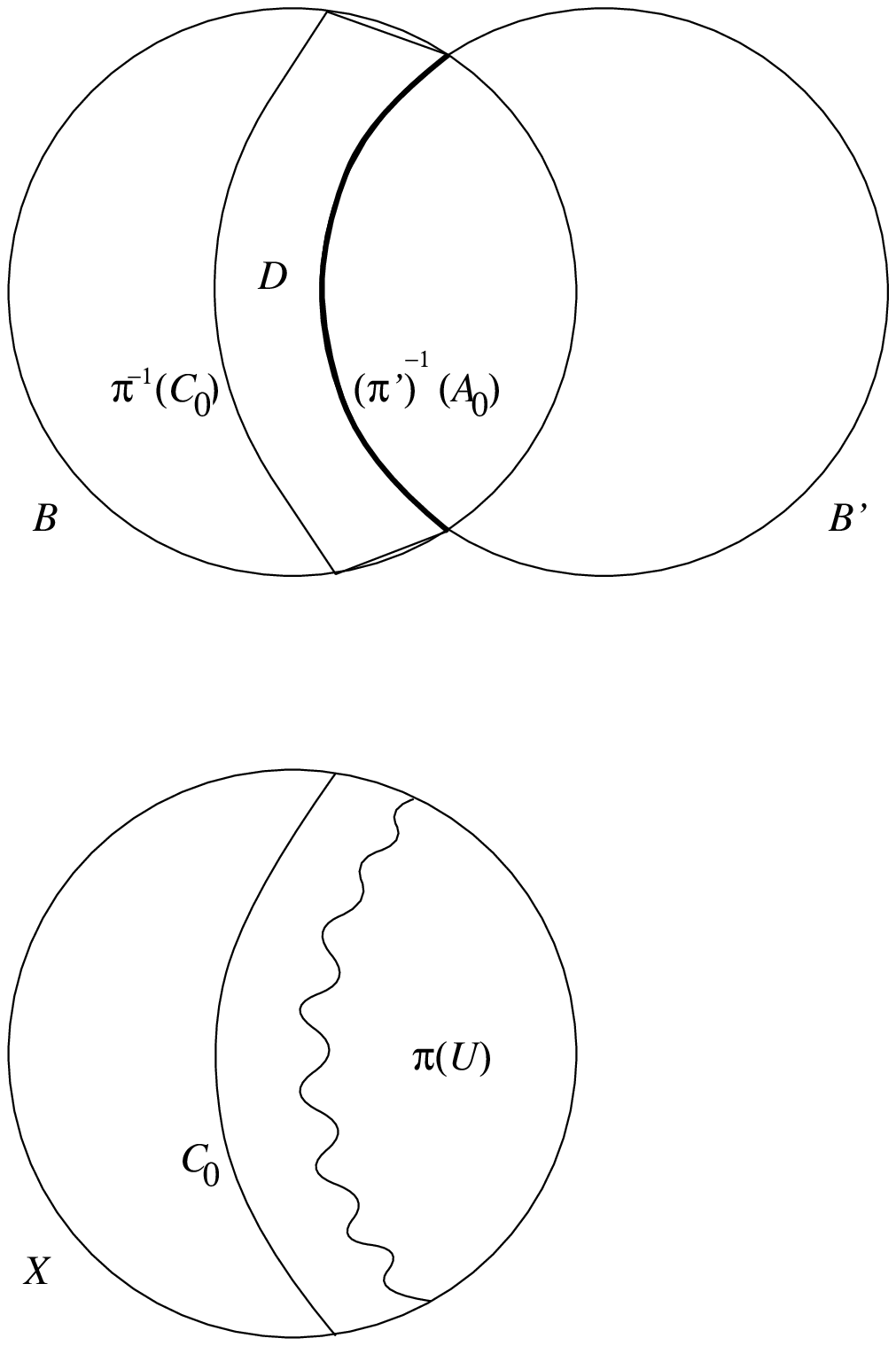}
  \caption{Case (B)}
  \label{fig:B}
\end{figure}
%
In this case $A_0$ is an arc.
In a way similar to Case (A), we apply Sublemma \ref{slem:isotopy}
to obtain a  closed domain $D$ bounded by $\pi^{-1}(C_0)$,
$(\pi')^{-1}(A_0)$ and two annuli joining  $(\pi')^{-1}(\partial A_0)$
and $(\pi')^{-1}(\partial C_0)$
which has an $S^1$-bundle structure via  
a bi-Lipschitz homeomorphism 
$k:I^2\times S^1\to D$ such that 
\begin{enumerate}
 \item for each $x\in \partial I^2$, $k(x\times S^1)$
      coincides with a fiber on $\partial D;$
 \item for each $x\in I^2$, $k(x\times S^1)$ is 
       $\nu$-perpendicular to horizontal subspaces.
\end{enumerate}
Thus we obtain the conclusion of Lemma \ref{lem:glue1}.
\end{proof}

Let $\nu>0$ be sufficiently small like $\nu=10^{-10}$, 
and let $\delta^*=\delta^*(\nu)$ 
and $\zeta=\zeta(\nu)$ be the constants given in Lemma \ref{lem:glue1}.
Letting $Q$ be the positive integer in Section \ref{sec:decomp} and 
setting
\begin{equation}
    \zeta^*:= \overbrace{\zeta(\zeta(\cdots (\zeta}^{%
                             \text{$Q$-times}}(\nu))\cdots)),\quad
          \delta^*:=\delta^*(\zeta^*),
                       \label{eq:delta}
\end{equation}
we choose $\epsilon$ in \eqref{eq:small-vol} satisfying
\begin{equation}
    \tau(2\delta^*) + \tau(s/2|\epsilon) < \zeta^*,
                       \label{eq:epsilon}
\end{equation}
where $r$  and $\sigma^*$ are defined as in \eqref{eq:r} and 
\eqref{eq:sigma}.
Note that $\zeta^*$, $\delta^*$, $\sigma^*$ and $s$ are universal
constants. We also choose a universal constant $\epsilon_1^*$ in 
Corollary \ref{cor:I} like $\epsilon_1^* \ll \sigma^*$.



\section{Proof of Lemma \ref{lem:glue}} \label{sec:glueII}

In this section, we shall prove Lemma \ref{lem:glue}.
$d_H^M$ denotes the Hausdorff distance in $M$.

Put 
\[
     B_{i_1,\ldots,i_k} := B_{i_1}\cup\cdots\cup B_{i_k},
\]
for $i_1,\ldots,i_k\in \cal I_2$.

\begin{ass} \label{ass:glue}
We assume that $B_{i_j}(1/10)$ meets $B_{i_1,\ldots,i_{j-1}}$ for 
every $2\le j\le k$. Then 
there exist a local $S^1$-action 
$\psi_{i_1,\ldots,i_{k}}$ on $\underline{B_{i_1,\ldots,i_{k}}}$
and a topological two-manifold $X_{i_1,\ldots,i_{k}}$
with the orbit projection 
$\pi_{i_1,\ldots,i_{k}}:\underline{B_{i_1,\ldots,i_{k}}} \to
X_{i_1,\ldots,i_{k}}$
satisfying the following $:$
\begin{enumerate}
 \item 
  \begin{equation*}
    \psi_{i_1,\ldots,i_k}= \begin{cases}
      \psi_{i_1,\ldots,i_{k-1}} &  \quad \text{on \,
              $\underline{B_{i_1}\cup\ldots\cup B_{i_{k-1}}-B_{i_k}}$}  \\
      \psi_{i_k} & \quad \text{on \,
              $\underline{B_{i_k}}$};
        \end{cases}
  \end{equation*}
\item $\underline{B_{i_1,\ldots,i_{k}}}$ is a manifold with boundary, and 
   $d_H^M(B_{i_1,\ldots,i_{k}},\underline{B_{i_1,\ldots,i_{k}}})
           < \tau(\epsilon)$.
 \item Each orbit of $\psi_{i_1,\ldots,i_{k}}$ has diameter
   $<\tau(\epsilon);$
 \item There are no singular orbits of $\psi_{i_1,\ldots,i_{k}}$
   over $\partial X_{i_1,\ldots,i_{k}} - \partial_{*} X_{i_1,\ldots,i_{k}}$,
   where $\partial_{*} X_{i_1,\ldots,i_{k}}$ is defined as in 
   \eqref{eq:bdy}.
 \item Let $N_{i_1,\ldots,i_k;1}\subset B_{i_1,\ldots,i_{k}}$ be the set of 
   points $p\in B_{i_1,\ldots,i_{k}}$ such that  
   there is a $(2, 2\delta^*)$-strainer at $p$ of length $\ge s/2$.
   Then every fiber $F$ of  $\psi_{i_1,\ldots,i_{k}}$ contained in 
   $N_{i_1,\ldots,i_k;1}$ is $\zeta^{Q-n}(\nu)$-perpendicular to 
   horizontal subspaces, where
    \begin{equation*}
       \zeta^{Q-n}(\nu) :=  
         \overbrace{\zeta(\cdots (\zeta}^{\text{$(Q-n)$-times}}(\nu))\cdots )),
    \end{equation*}
   and $n$ denotes the number of $\hat B_i$ containing $F$.
\end{enumerate}
\end{ass}

\begin{proof}
We prove it by induction. Assertion \ref{ass:glue} certainly 
holds for $k=1$ by Theorem \ref{thm:loc-fib}, 
Lemmas \ref{lem:fiber}(2) and \ref{lem:N_1}.
Assume that a local
$S^1$-action $\psi_{i_1,\ldots,i_{k-1}}$ on 
$\underline{B_{i_1,\ldots,i_{k-1}}}$  satisfying Assertion \ref{ass:glue}
has been constructed.
For simplicity, we set
\begin{gather*}
  B:=\underline{B_{i_1,\ldots,i_{k-1}}}, \quad X:=X_{i_1,\ldots,i_{k-1}},\\
  \hat B':= \hat B_{i_k},\quad B':=B_{i_k},
       \quad \hat X':= \hat X_{i_k}, \quad X':= X_{i_k}, \\
   \psi:=\psi_{i_1,\ldots,i_{k-1}}, \quad \psi':=\psi_{i_k} \\
    \pi:=\pi_{i_1,\ldots,i_{k-1}}, \quad \pi':=\pi_{i_k}.
\end{gather*}
We shall carry out a gluing procedure on each component, say $U$,
of $B\cap B'$ using Lemma \ref{lem:glue1}.
$\partial_{*}X$ and $\partial_{*}X'$ are defined as in 
\eqref{eq:bdy}.
Let $A\subset \partial X'$, $A_0$, $C\subset X$ and $C_0$
be defined as in the proof of Lemma
\ref{lem:glue1}.

First consider the case when $A_0$ does not meet 
$B(S_{\delta^*}(Z'),r)$, where $Z'$ is the Alexandrov space
containing $X'$, which implies  $(\pi')^{-1}(A_0)\subset N'$.
Thus every $\pi'$-fiber in  $(\pi')^{-1}(A_0)$ is 
$\zeta^{Q}(\nu)$-perpendicular to horizontal subspaces.
Since $\pi^{-1}(C_0)$ is close to $(\pi')^{-1}(A_0)$, we obtain that 
$\pi^{-1}(C_0)\subset N_{i_1,\ldots,i_{k-1};1}$. 
Condition (5) of Assertion \ref{ass:glue} for
$\psi_{i_1,\ldots,i_{k-1}}$ then implies that 
every $\pi$-fiber in  $(\pi)^{-1}(C_0)$ is 
$\zeta^{Q-n+1}(\nu)$-perpendicular to horizontal subspaces.
Therefore we can apply Lemma \ref{lem:glue1} to get the 
required gluing of $\psi$ and $\psi'$ on a neighborhood 
joining $(\pi')^{-1}(A_0)$ and $\pi^{-1}(C_0)$.

Next consider the other case when  $A_0$  meets
$B(S_{\delta^*}(Z'), r)$. In this case 
it follows from the condition (c) of Case (B) of Theorem \ref{thm:loc-fib} 
that an endpoint of $A_0$ must be contained in 
$B(\partial_*X', r)$.

\begin{figure}[htbp]
  \centering
  \includegraphics[scale=0.8]{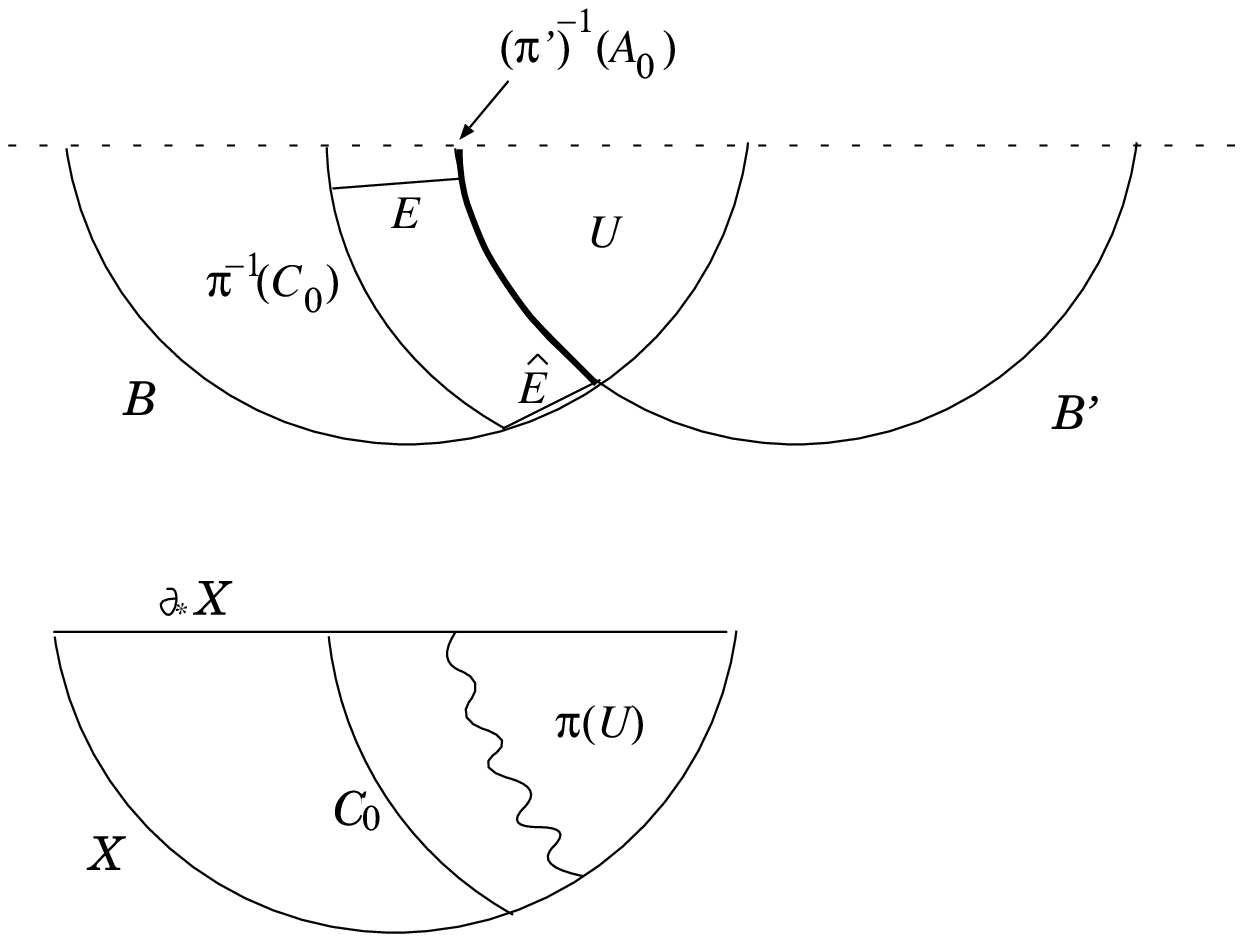}
  \caption{}
  \label{fig:C}
\end{figure}

%
%
Take $x, \hat x\in A_0$ such that the subarc $[x,\hat x]$ of $A_0$
is contained in $X'-B(S_{\delta^*}(Z'),r)$ and $A_0-[x,\hat x]$ 
as well as $x$ is contained in $B(\partial_*X',r)$.
The case $\hat x\in B(\partial_*X',r)$ may happen. 
Let $y, \hat y\in C_0$ be defined as in the proof of Lemma \ref{lem:glue1}
for $x$, $\hat x$.
Applying (the proof of) Lemma \ref{lem:glue1}, we have a 
domain $D$ bounded by $(\pi')^{-1}([x,\hat x])$, 
$\pi^{-1}([y,\hat y])$ and two annuli, say $E$ and $\hat E$, 
having a compatible 
$S^1$-fiber structure via a bi-Lipschitz homeomorphism  
$k:D^2\times S^1\to D$ such that 
for each $x\in D^2$, $k(x\times S^1)$ is 
$\zeta^{Q-n}(\nu)$-perpendicular to horizontal subspaces.
Note that each component, say $W$,  of the domain bounded by 
$(\pi')^{-1}(A_0-[x,\hat x])$, $\pi^{-1}(C_0-[y,\hat y])$ and 
$E$ (and possibly $\hat E$) is a three-disk. Therefore
one can put a compatible structure of $S^1$-action on $W$
all of whose orbit has diameter $<\tau(\epsilon)$.

From the gluing constructions above,
we obtain 
$\underline{B_{i_1,\ldots,i_{k}}}$ and a local $S^1$-action 
$\psi_{i_1,\ldots,i_{k}}$ on it satisfying the conclusion of 
the assertion.
This completes the proof of Assertion \ref{ass:glue}.
\end{proof}

From \eqref{eq:biggest} and the gluing argument 
used in the proof of Assertion \ref{ass:glue}, it is now obvious
that $\underline{U_2}\simeq  U_2$.
Thus we have completed the proof of Lemma \ref{lem:glue}.



\section{Thick-thin decomposition} \label{sec:thick}

For a three-manifold $N$, we denote by $\CCap(N)$ a three-manifold 
obtained by gluing of $N$ and some copies of $D^3$ along
all the sphere-components of $\partial N$.

\begin{thm} \label{thm:thick}
If $M$ is a closed orientable Riemannian three-manifold with
sectional curvature $K\ge -1$, then we have a decomposition
\[
    M = M_{\rm thick} \cup M_{\rm thin},
\]
satisfying the following$:$
Let $\epsilon_0$ and $\delta_0$ be positive numbers given 
Theorem \ref{thm:graph}.
 \begin{enumerate}
  \item For every $p\in M_{\rm thick}$, $\vol(B(p,1))\ge \epsilon_0/2$.
  \item $\CCap(M_{\rm thin})$ is homeomorphic to a graph manifold
    if $\diam(M)\ge\delta_0$.
 \end{enumerate}
\end{thm}

Roughly speaking, $M_{\rm thin}$ is a piece of $M$ which 
collapses. 

Theorem \ref{thm:thick} is closely related with a result in 
\cite{Pr:surgery}, where a thick-thin 
decomposition of a closed three-manifold in connection with Ricci flow 
is announced.
\par

\begin{proof}[Proof of Theorem \ref{thm:thick}]
If $\diam(M)<\delta_0$, then we put $M_{\rm thin}:=M$.
For the proof of Theorem \ref{thm:thick},
we may assume $\diam(M)\ge\delta_0$.
For every point $p\in M$, one of the following holds:
\begin{enumerate}
 \item $\vol B(p,1) \ge \epsilon_0$. In this case,
  we put $B_p:= B(p,1)$ and  $X_p:= B(p,1)$.
 \item $\vol B(p,1) < \epsilon_0$.  In this case, 
  by Theorem \ref{thm:loc-fib} there exist a small perturbation $B_p$ of 
  $B(p,1)$ and a metric ball $X_p$ in some complete Alexandrov space $X$ 
  with curvature $\ge -1$ and $1\le \dim X\le 2$
  such that $B_p$ has a fiber structure over $X_p$.
\end{enumerate}
Take points $p_1, p_2, \ldots,$ of $M$ as in Section \ref{sec:decomp}
such that the collection $\{ B_{p_i}\}$ given as above is a 
covering of $M$.
Let $X_{p_i}$ be also chosen as above.
For simplicity, we put 
 \[ 
       B_i:=B_{p_i},\quad  X_i:=X_{p_i}.
 \]
If $\dim X_i=2$, then there exists a local 
$S^1$-action $\psi_i$ on $B_i$ such that 
$B_i/\psi_i \simeq X_i$.
For each $1\le j\le 3$,
let $\cal I_j$ denote the set of all $i$ with $\dim X_i=j$, 
and consider 
\[
    U_j := \bigcup_{i\in \cal I_j} B_i.
\]
By Lemm \ref{lem:dim1}, each component of $U_1$ is either 
cylindrical or cylindrical with a cap.
By Lemma \ref{lem:glue}, we have a small perturbation
$\underline{U_2}$ of $U_2$ on which one can construct a local $S^1$-action. 
Note that each component of $M_{\rm thin}:= U_1 \cup \underline{U_2}$
is homeomorphic to $S^2$ or $T^2$. The argument in Section \ref{sec:decomp}
shows that $\CCap(M_{\rm thin})$ is homeomorphic to a graph manifold.
Now it is obvious that 
for any point $p$ in $M_{\rm thick}:= \overline{M- M_{\rm thin}}$,
$\vol(B(p,1))\ge \epsilon_0/2$.
\end{proof}

\section{Appendix: Collapsing under a local lower curvature bound} 
       \label{sec:loc-col}

In this appendix, we give a short description about collapsing
three-manifold under a local lower curvature bound, which is 
discussed in \cite{Pr:surgery}.

For a positive number $\epsilon$, a closed Riemannian $n$-manifold 
$(M,g)$ is called {\em $\epsilon$-collapsed under a local lower 
curvature bound} 
if for each  $x \in M$, there exists a 
$\rho$, $0 < \rho \le \diam(M,g)$,  with 
\begin{equation}
     \vol B(x,\rho) \le \epsilon\rho^n, \qquad
          K \ge -\rho^{-2} \quad  \text{on $B(x,\rho)$}.
               \label{eq:loc-col}
\end{equation}

Theorem \ref{thm:graph} extends to the following:

\begin{thm}[Theorem 7.4 \cite{Pr:surgery}] \label{thm:loc-col}
Let $\epsilon_0$ be a positive number given in Theorem \ref{thm:graph}.
If a closed orientable Riemannian three-manifold $(M,g)$ is 
$\epsilon_0$-collapsed 
under a local lower curvature bound, then it is homeomorphic to 
a graph manifold.
\end{thm}

Let $\rho(x)$,  $0 < \rho(x) \le \diam(M,g)$, be the supremum 
of $\rho>0$ satisfying \eqref{eq:loc-col}.
By Theorem \ref{thm:loc-fib},  a small perturbation $B_x$ of
$B(x,\rho(x))$ has a singular fibration over a metric ball in 
some Alexandrov space
with curvature $\ge -1$ and dimension one or two.
Choose a covering $\{ B_{x_i} \}$ of $M$
such that $\{ B(x_i, \rho(x_i)/10)\}$ is a maximal disjoint family.
Let $U_1$, $U_2$ and $M=U_1\cup \hat U_2$ be as in Section \ref{sec:decomp}.
In a way similar to Lemma \ref{lem:dim1}, one can prove that each component
of $U_1$ is either cylindrical or cylindrical with a cap.

Let $\mathcal B_2$ be defined as in Section \ref{sec:decomp}.

\begin{lem} \label{lem:rho}
Suppose $B_x$ and $B_y$ in $\mathcal B_2$
satisfy that $B(x,\rho(x)/4)$ meets $B(y,\rho(y)/4)$.
Then 
\[
       C^{-1} \le \frac{\rho(x)}{\rho(y)} \le C,
\]
for some universal positive number $C$.
\end{lem}

\begin{proof} 
Assuming $\rho(x)< \rho(y)$, we put $\rho(y) = R \rho(x)$.
By triangle inequality, 
$B(x, R_1\rho(x))\subset B(y,\rho(y))$,
where $R_1:=R/2$, which implies that 
\begin{equation}
    K \ge - (R_1\rho(x))^{-2}\quad \text{on\, $B(x,R_1\rho(x))$}.
       \label{eq:curvature}
\end{equation}
Next we show 
\begin{equation}
   \vol B(x, R_1\rho(x)) \le \epsilon_0 (R_1\rho(x))^3
          \label{eq:large-ball}
\end{equation}
if $R_1$ is larger than some uniform positive constatn.
In view of the maximality of $\rho(x)$, \eqref{eq:curvature} and 
\eqref{eq:large-ball} yield the conclusion of the lemma.

First note that large parts of $B_x$ and $B_y$ have
$S^1$-fiber structures.
Let $\ell$ denote the length of a regular circle fiber $F$ 
contained in $B(x,\rho(x)/4)\cap B(y,\rho(y)/4)$.
Let $g_x:=\rho(x)^{-2} g$.
Since $B(y,\rho(y)/8)\subset B(x,R_1\rho(x))\subset B(y,\rho(y))$,
Lemma \ref{lem:length} implies 
\[
   C_1^{-1}\frac{\ell}{\rho(y)} \le \vol_{g_y} B(x,R_1\rho(x))\le
                           C_1\frac{\ell}{\rho(y)},
\]
for some uniform positive number $C_1$. It follows that 
\[
 4C_1^{-1}\ell R_1^2 \rho(x)^2 \le \vol_{g}B(x,R_1\rho(x))
           \le 4C_1\ell R_1^2 \rho(x)^2.
\]
Hence \eqref{eq:large-ball} holds if $4C_1\ell < \epsilon_0 R_1\rho(x)$.
On the other hand, it follows from the assumption and Lemma \ref{lem:length}
that 
\[
 C_1^{-1}\frac{\ell}{\rho(x)} \le \vol_{g_x} B(x,\rho(x))< \epsilon_0.
\]
Therefore we obtain \eqref{eq:large-ball} for $R_1>4C_1^2$.

We conclude that the constant $C$ in the lemma is given 
by $C=8C_1^2$.
\end{proof}

Lemma \ref{lem:rho} together with
the Bishop-Gromov comparison theorem yields a uniform 
upper bound on the maximal number of intersections among
the metric balls in $\mathcal B_2$. Therefore 
our local gluing argument in Sections \ref{sec:glueI} and  \ref{sec:glueII}
goes through the present context as well 
to complete the proof of Theorem \ref{thm:loc-col}. 



\bibliographystyle{amsplain}

\bibliography{all}

\providecommand{\bysame}{\leavevmode\hbox to3em{\hrulefill}\thinspace}
\providecommand{\MR}{\relax\ifhmode\unskip\space\fi MR }
\providecommand{\MRhref}[2]{%
  \href{http://www.ams.org/mathscinet-getitem?mr=#1}{#2}
}
\providecommand{\href}[2]{#2}
\begin{thebibliography}{10}

\bibitem{BGP}
Yu. Burago, M.~Gromov, and G.~Perel'man, \emph{A. {D}. {A}leksandrov spaces
  with curvatures bounded below}, Uspekhi Mat. Nauk \textbf{47} (1992),
  no.~2(284), 3--51, 222, translation in Russian Math. Surveys 47 (1992), no.
  2, 1--58.

\bibitem{CG:collapseI}
J.~Cheeger and M.~Gromov, \emph{Collapsing {R}iemannian manifolds while keeping
  their curvature bounded {I}}, J. Differential Geom. \textbf{23} (1986),
  no.~3, 309--346.

\bibitem{CG:collapseII}
\bysame, \emph{Collapsing {R}iemannian manifolds while keeping their curvature
  bounded {II}}, J. Differential Geom. \textbf{32} (1990), no.~1, 269--298.

\bibitem{FY:fundgp}
K.~Fukaya and T.~Yamaguchi, \emph{The fundamental groups of almost
  nonnegatively curved manifolds}, Ann. of Math. (2) \textbf{136} (1992),
  253--333.

\bibitem{Or:Seif}
P.~Orlik, \emph{Seifert manifolds}, Lecture Notes in Math., no. 291,
  Springer-Verlag, Berlin-New York, 1972.

\bibitem{Pr:alex2}
G.~Perelman, \emph{{A}. {D}. {A}lexandrov's spaces with curvatures bounded from
  below {II}}, preprint.

\bibitem{Pr:surgery}
\bysame, \emph{Ricci flow with surgery on three-manifolds}, arXiv:math. DG /
  0303109.

\bibitem{ST:cut}
K.~Shiohama and M.~Tanaka, \emph{Cut loci and distance spheres on {A}lexandrov
  surfaces}, Actes de la Table Ronde de G\'eom\'etrie Diff\'erentielle (Luminy,
  1992), S\'emin. Congr., vol.~1, Soc. Math. France, Paris, 1996, pp.~531--559.

\bibitem{SyYm:3mfd}
T.~Shioya and T.~Yamaguchi, \emph{Collapsing three-manifolds under a lower
  curvature bound}, J. Differential Geom. \textbf{56} (2000), no.~1, 1--66.

\bibitem{Ym:collapsing}
T.~Yamaguchi, \emph{Collapsing and pinching under a lower curvature bound},
  Ann. of Math. (2) \textbf{133} (1991), 317--357.

\bibitem{Ym:conv}
\bysame, \emph{A convergence theorem in the geometry of {A}lexandrov spaces},
  Actes de la Table Ronde de G\'eom\'etrie Diff\'erentielle (Luminy, 1992),
  S\'emin. Congr., vol.~1, Soc. Math. France, Paris, 1996, pp.~601--642.

\end{thebibliography}

\end{document}